
\def\LaTeX{%
  \let\Begin\begin
  \let\End\end
  \let\salta\relax
  \let\finqui\relax
  \let\futuro\relax}

\def\UK{\def\our{our}\let\sz s}
\def\USA{\def\our{or}\let\sz z}

\UK



\LaTeX

\USA


\salta

\documentclass[twoside,12pt]{article}
\setlength{\textheight}{24cm}
\setlength{\textwidth}{16cm}
\setlength{\oddsidemargin}{2mm}
\setlength{\evensidemargin}{2mm}
\setlength{\topmargin}{-15mm}
\parskip2mm


\usepackage[usenames,dvipsnames]{color}
\usepackage{amsmath}
\usepackage{amsthm}
\usepackage{amssymb}
\usepackage[mathcal]{euscript}
\usepackage{enumitem}
%
%
\usepackage{cite}
%
%
%
%


\definecolor{viola}{rgb}{0.3,0,0.7}
\definecolor{ciclamino}{rgb}{0.5,0,0.5}
\definecolor{rosso}{rgb}{0.8,0,0}

\def\dafare #1{{\color{red}#1}}

\def\dafare #1{#1}


\bibliographystyle{plain}


%

\finqui

\def\Beq{\Begin{equation}}
\def\Eeq{\End{equation}}
\def\Bsist{\Begin{eqnarray}}
\def\Esist{\End{eqnarray}}
\def\Bdef{\Begin{definition}}
\def\Edef{\End{definition}}
\def\Bthm{\Begin{theorem}}
\def\Ethm{\End{theorem}}
\def\Blem{\Begin{lemma}}
\def\Elem{\End{lemma}}

\def\Bcor{\Begin{corollary}}
\def\Ecor{\End{corollary}}
\def\Brem{\Begin{remark}\rm}
\def\Erem{\End{remark}}

\def\Bdim{\Begin{proof}}
\def\Edim{\End{proof}}
\def\Bcenter{\Begin{center}}
\def\Ecenter{\End{center}}
\let\non\nonumber




\def\step #1 \par{\medskip\noindent{\bf #1.}\quad}


\def\Lip{Lip\-schitz}
\def\Holder{H\"older}
\def\Frechet{Fr\'echet}

\def\aand{\quad\hbox{and}\quad}

\def\wk{well-known}
\def\socal{so-called}
\def\lhs{left-hand side}
\def\rhs{right-hand side}
\def\sfw{straightforward}

\def\CH{Cahn-Hilliard}



\def\multibold #1{\def\arg{#1}%
  \ifx\arg\pto \let\next\relax
  \else
  \def\next{\expandafter
    \def\csname #1#1#1\endcsname{{\bf #1}}%
    \multibold}%
  \fi \next}

\def\pto{.}

\def\multical #1{\def\arg{#1}%
  \ifx\arg\pto \let\next\relax
  \else
  \def\next{\expandafter
    \def\csname cal#1\endcsname{{\cal #1}}%
    \multical}%
  \fi \next}


\def\multimathop #1 {\def\arg{#1}%
  \ifx\arg\pto \let\next\relax
  \else
  \def\next{\expandafter
    \def\csname #1\endcsname{\mathop{\rm #1}\nolimits}%
    \multimathop}%
  \fi \next}

\multibold
qwertyuiopasdfghjklzxcvbnmQWERTYUIOPASDFGHJKLZXCVBNM.

\multical
QWERTYUIOPASDFGHJKLZXCVBNM.

\multimathop
ad dist div dom meas sign supp .


\def\accorpa #1#2{\eqref{#1}--\eqref{#2}}
\def\Accorpa #1#2 #3 {\gdef #1{\eqref{#2}--\eqref{#3}}%
  \wlog{}\wlog{\string #1 -> #2 - #3}\wlog{}}


\def\graffe #1{\mathopen\{#1\mathclose\}}

\def\<#1>{\mathopen\langle #1\mathclose\rangle}
\def\norma #1{\mathopen \| #1\mathclose \|}

\def\intQtT{\int_{Q_t^T}}
\def\intQ{\int_Q}
\def\iO{\int_\Omega}

\def\inttt{\int_t^T}

\def\dt{\partial_t}
\def\dn{\partial_n}

\def\checkmmode #1{\relax\ifmmode\hbox{#1}\else{#1}\fi}

\def\aeQ{\checkmmode{a.e.\ in~$Q$}}

\def\aaQ{\checkmmode{for a.a.~$(x,t)\in Q$}}

\def\aat{\checkmmode{for a.a.~$t\in(0,T)$}}


\def\erre{{\mathbb{R}}}

\def\enne{{\mathbb{N}}}




\def\genspazio #1#2#3#4#5{#1^{#2}(#5,#4;#3)}
\def\spazio #1#2#3{\genspazio {#1}{#2}{#3}T0}

\def\L {\spazio L}
\def\H {\spazio H}

\def\C #1#2{C^{#1}([0,T];#2)}


\def\Lx #1{L^{#1}(\Omega)}
\def\Hx #1{H^{#1}(\Omega)}

\def\Ldue{\Lx 2}

\def\Huno{\Hx 1}
\def\Hdue{\Hx 2}



\let\theta\vartheta

\let\phi\varphi

\let\TeXchi\chi                         
\newbox\chibox
\setbox0 \hbox{\mathsurround0pt $\TeXchi$}
\setbox\chibox \hbox{\raise\dp0 \box 0 }
\def\chi{\copy\chibox}



\def\bQ{b_1}
\def\bO{b_2}
\def\bQh{b_3}
\def\bOh{b_4}
\def\bz{b_0}

\def\phQ{\phi_Q}

\def\phO{\phi_\Omega}


\def\sO{\s_\Omega}

\def\sQ{\s_Q}

\def\Uad{\calU_{\ad}}
\def\uopt{\overline u}

\def\Vp{V^*}

\let\hat\widehat

\def\Pi{\hat\pi}


\def\cd{c_\delta}
\def\s{\sigma}  
\def\m{\mu}	    
\def\ph{\phi}	
\def\a{\alpha}	
\def\b{\beta}	
\def\d{\delta}  
\def\bph{\overline\ph}  
\def\bm{\overline\m}    
\def\bs{\overline\s}    
\def\J{{\cal J}}
\def\S{{\cal S}}
\def\I2 #1{\int_{Q_t}|{#1}|^2}
\def\IN2 #1{\int_{Q_t}|\nabla{#1}|^2}
\def\IO2 #1{\iO |{#1(t)}|^2}
\def\INO2 #1{\iO |\nabla{#1}(t)|^2}
\def\UR{{\cal U}_R}
\def\CP{\boldsymbol{(CP)}}



\Begin{document}

\title{
		Vanishing parameter for an optimal control problem
		\\[0.3cm] 
		modeling tumor growth}
\author{}
\date{}
\maketitle

\Bcenter
\vskip-1cm
{\large\sc Andrea Signori$^{(1)}$}\\
{\normalsize e-mail: {\tt andrea.signori02@universitadipavia.it}}\\[.25cm]
$^{(1)}$
{\small Dipartimento di Matematica e Applicazioni, Universit\`a di Milano--Bicocca}\\
{\small via Cozzi 55, 20125 Milano, Italy}

\Ecenter
\Begin{abstract}\noindent
A distributed optimal control problem for a phase field
system which physical context is that of tumor growth is discussed.
The system we are going to take into account
consists of a Cahn-Hilliard equation for the phase variable
(relative concentration of the tumor) coupled
with a reaction-diffusion equation for the nutrient.
The cost functional is of standard tracking-type and the control variable models 
the intensity at which it is possible to dispense medication.
The model we deal with presents two small and positive parameters which
are introduced in previous contributions as relaxation terms.
Here, starting from the already investigated optimal control problem
for the relaxed model,
we aim at confirming the existence of optimal control and
characterizing the first-order necessary optimality condition, via asymptotic schemes,
when one of the two occurring parameters goes to zero.
%

\vskip3mm
\noindent {\bf Key words}
Asymptotic analysis, distributed optimal control, phase field model, 
tumor growth, cancer treatment, evolution equations, Cahn-Hilliard equation, 
optimal control, adjoint system, necessary optimality conditions.
\vskip3mm
\noindent {\bf AMS (MOS) Subject Classification} 
35K61,  
35Q92,  
49J20,  
49K20,  
35K86,  
92C50.  
\End{abstract}

\vskip3mm

\pagestyle{myheadings}
\newcommand\testopari{\sc {Vanishing parameter for an optimal control problem}}
\newcommand\testodispari{\sc Signori}
\markboth{\testodispari}{\testopari}

\salta
\finqui

\newpage
\section{Introduction}
\label{SEC_INTRODUCTION}
\setcounter{equation}{0}
The need for investigating tumor growth from a mathematical viewpoint
stems from the great impact that it may have on medical treatments.
As a matter of fact, in the last years, there is increasing attention
by the mathematical community toward biological and medical models (see, e.g., \cite{CL}). 
In particular, an open and unknown area such as the tumor field
can find a useful support tool in the mathematical predictions.
In fact, this latter 
could be able to pull out some of the main features of the evolution
phenomena and, by focusing on some particular aspects, it may give some deep
insights as if a given negative outcome was to be foreseen, it would be possible to prevent it.
Moreover, the theoretical investigation has the
huge advantage that no patient is put at risk.
Furthermore, without the claim to cure the disease, the mathematical models 
could provide prominent a priori information as a support for the medical treatments
leading to more personalized therapy. 
Indeed, despite the wide number of parameters involved in the disease,
due to the few understanding 
of the tumor evolution, the corresponding clinical treatment is quite standardized,
while every patient responds differently to the medications.

Among the numerous models recently proposed, we focus on the
ones derived by continuum mixture and phase field theories.
The evolution of a young tumor, before the development of quiescent cells,
can be described as a \CH~equation for the phase variable
(see, e.g., \cite{Mir_CH} and the huge references therein for a general, while rich,
introduction to the \CH~equation), 
coupled with a reaction-diffusion for an unknown
species acting as a nutrient (e.g., oxygen or glucose).
The model we are going to face in this work consists of a variation of the one 
introduced by Hawkins-Daruud et al. in \cite{HDZO}, where
the velocity contributions are neglected (see also \cite{HKNZ,WZZ,CLLW,HDPZO,OHP}).
Several models, by interpreting the tumors and the healthy cells as inertia-less fluids,
also include the contribution of the velocity field assuming
a Darcy law or a Stokes-Brinkman equation. In this regards, let us refer to
\cite{WLFC,GLSS,DFRGM ,GARL_1,GARL_4,GAR, EGAR, FLRS,GARL_2, GARL_3}, 
where further mechanisms such as active transport and chemotaxis 
are also taken into account. 
We also point out the paper \cite{FRL}, where a non-local model is proposed.

At first, let us point out that
the symbol $\Omega \subset \erre^3$ is devoted to indicating the 
set where the evolution takes place which boundary we denote by $\Gamma$.
Furthermore, given a final time $T>0$, we set for convenience
\Bsist
	& \non
	Q_{t}:=\Omega \times (0,t), \quad \Sigma_{t}:=\Gamma \times (0,t)
	\quad \hbox{for every $t \in (0,T]$,}
	\\
	& Q:=Q_{T}, \quad \hbox{and} \quad \Sigma:=\Sigma_{T}.
	\label{QS}
\Esist

In the present paper, we are going to deal with the optimal control problem 
consisting of minimizing
the \socal~objective, or tacking-type, cost functional  
\begin{align}
	\non
	\J (\ph, \s, u)   &:= 
	\frac \bQ 2 \norma{\ph - \phQ}_{L^2(Q)}^2
	+\frac \bO 2 \norma{\ph(T)-\phO}_{L^2(\Omega)}^2
	+ \frac \bQh 2 \norma{\s - \sQ}_{L^2(Q)}^2
	\\ & \qquad
	+ \frac \bOh 2 \norma{\s(T)-\sO}_{L^2(\Omega)}^2
    + \frac \bz 2 \norma{u}_{L^2(Q)}^2,
    \label{costfunct}
\end{align}
subject to the control-box constraints
\Beq
    \label{Uad}
	u \in \Uad := \graffe{u \in L^{\infty}(Q): u_* \leq u \leq u^*\ \aeQ},
\Eeq
and under the assumption that the variables $\ph$ and $\s$
solve the following system
\begin{align}
   \a \dt \m + \dt \ph - \Delta \m &= P(\phi) (\sigma - \mu)
  \quad \hbox{in $\, Q$}
  \label{EQold_prima}
  \\
   \mu &= \beta \dt \phi - \Delta \phi + F'(\phi)
  \label{EQold_seconda}
  \quad \hbox{in $\,Q$}
  \\
   \dt \sigma - \Delta \sigma &= - P(\phi) (\sigma - \mu) + u
  \label{EQold_terza}
  \quad \hbox{in $\,Q$}
  \\
   \dn \m &= \dn \ph =\dn \s = 0
  \quad \hbox{on $\,\Sigma$}
  \label{EQold_BC}
  \\
   \m(0)&=\m_0,\, \ph(0)=\ph_0,\, \s(0)=\s_0
  \quad \hbox{in $\,\Omega.$}
  \label{EQold_IC}
\end{align}
\Accorpa\EQold EQold_prima EQold_IC
For the sake of synthesis, let us describe the physical
background of the occurring variables without diving into the details.
%
The admissible set $\Uad$ fix the space in which the control variable $u$ can be chosen
and it is given in terms of the bounds $u_*$ and $u^*$.
Moreover, $\bz, \bQ, \bO, \bQh,\bOh$ stand for nonnegative constants, not all zero,
while $\phQ, \sQ, \phO, \sO $ denote some target functions defined in $Q$ and 
$\Omega$, respectively.
The variable $\ph$ is an order parameter and it is designed to keep track of the 
evolution of the tumor in the tissue. It is a normalized relative concentration
and it ranges between $-1$ and $+1$, where 
these extremes represent the pure phases, that is the tumorous and the
healthy case, respectively.
Furthermore, the variable $\m$ stands, as usual for the \CH~equation,
for the chemical potential for $\ph$. The third unknown $\s$ has the role of
describing the evolution of the nutrient within the evolution process and
it is normalized between $0$ and $1$ with the following property:
the closer to one, the richer of nutrient the extra-cellular is, 
while the closer to zero, the poorer it is. Lastly, the variable $u$ represents
the \socal~control variable and, since it appears in the nutrient equation, 
it can be read as a supply of a nutrient or a drug in the medical treatment.
As the functions $P$ and $F$ are concerned, they are nonlinearities. The former models 
the proliferation of the tumor, while the latter
is a double-well potential associated with the \CH~equation.
A typical example of $F$ is the regular potential which is defined as follows
\Bsist
  & F_{reg}(r) = \frac14(r^2-1)^2 
  = \frac 14 ((r^2 - 1)^+)^2 + \frac 14 ((1 - r^2)^+)^2 \, 
  \quad \hbox{for } r \in \erre.
  \label{F_reg}
\Esist
We will see that the optimal control problem we are going to deal with 
will demand to restrict
the analysis on potentials which slightly generalize \eqref{F_reg}.
Namely, we cannot take into account singular or non-regular potential as the \wk~logarithmic
double-well potential or the double-obstacle one.
For different physically meaningful choices of the potentials, 
we refer to \cite{Agosti} and
to the references therein, where numerical simulations and comparison with
clinical data can be found as well.
Further details regarding the interpretation of the 
model can be found in \cite{FGR, CGH,CGRS_VAN,CGRS_ASY}.

The above system has already been investigated 
in \cite{CGH}, in the case where $\a=\b>0$, from the viewpoint of well-posedness and
long-time behavior in terms of the omega-limit set.
Further comprehension of the model has been achieved by \cite{CGRS_VAN,CGRS_ASY},
where the authors show under which framework the parameters $\a$ and $\b$
can be let to zero and also point out the existence and
uniqueness of the solution to the limit problem in their natural setting.
In addition, we also refer to \cite{FGR} where the system formally obtained by imposing 
$\a=\b=0$ is tackled. There, after providing the well-posedness, the authors 
focus on the long-time behavior of the solution in terms of the global attractor
(see, e.g., \cite{MirZel} for details on the asymptotic behavior of 
infinite-dimensional dynamical systems).
As for the long-time behavior of the same system, namely \EQold~with
$\a=\b=0$, we are also aware of the recent contributions \cite{MRS,CRW}.
Lastly, let us mention \cite{Kur}, where the author confirms the existence
of the above problem \EQold~when $\b\searrow 0$, extending the analysis to
the case of unbounded domains by accounting for suitable approximation schemes.

As the terminology is concerned, the system that the control variable has to satisfy
is referred to as the state system. 
Moreover, once that the  well-posedness of the state system has been performed, 
we can introduce the control-to-state
mapping as the map that assigns to a given control the associated solution, 
namely the function
\Beq
	\non
	\S : u \mapsto \S(u):=(\m,\ph,\s)=(\m(u),\ph(u),\s(u)).
\Eeq
Starting from this, one can interpret the cost functional $\J$ as a function
depending on the control variable only, giving rise to the \socal~reduced cost functional
reading as
\Beq
	\non
	\J_{red} (u):= \J(\S_2(u),\S_3(u),u), 
\Eeq
where $\S_2$ and $\S_3$ denote the second and third component of the 
solution operator $\S$, respectively.

Even though the literature around the mathematical investigations of
biological and medical models find several examples, the corresponding 
optimal control contributions are very few.
Up to our knowledge, the first paper dealing with an optimal control problem for
a system very close the one gave above, namely the case $\a=\b=0$, is \cite{CGRS_OPT}.
Furthermore, we mention \cite{S}, which is our starting point.
There, the author handles the optimal control
problem for the classical tracking-type cost functional in the non-trivial case of the
logarithmic potential, where the presence of the relaxation terms turns out to be fundamental.
Moreover, in a following work, the same author proves that, 
accounting for an asymptotic technique 
known in the literature as to deep quench limit, it is also possible to generalize
the assumptions for the potentials in order to take into account also singular and 
non-regular potentials like the double-obstacle one. 
In addition, we refer to \cite{S_b}, where it was shown that the optimal control
problem for the state system \EQold~with $\b=0$ can be solved, by 
letting $\b\searrow 0$ in the optimal control problem associated with \EQold.
Lastly, let us address to \cite{CRW}, where an optimal control problem for 
\EQold, with $\a=\b=0$, has been discussed for a slightly more general class of cost functional
which takes into account the time optimization (see also
\cite{GARLR}, where the same generalized cost functional is taken for a different 
state system). 
To conclude the overview concerning the literature, let us
also point out \cite{SM}, where a different kind of control problem, known as sliding
mode control, is performed.
As for different state systems, let us refer to the recent
\cite{EK, EK_ADV}, where the authors establish the
existence of optimal controls and also characterize the optimality conditions
for the more involved \CH-Brinkman equation. Lastly, let us mention \cite{SW},
where a distributed optimal control problem for
the Cahn-Hilliard-Darcy system with mass source was studied.

Here, we aim to employ an asymptotic scheme similar to the one of \cite{S_b},
by letting $\a\searrow 0$ instead of $\b$, and assuming \cite{S}. 
Note that the present contribution complete the
picture around the optimal control problem for system \EQold~with
the standard tracking-type cost functional. Indeed,
the case $\a,\b>0$ has been investigated in \cite{S}, the case 
$\a>0$ and $\b=0$ has been studied in \cite{S_b}, whereas the case $\a=\b=0$ has been
treated in \cite{CGRS_OPT}.

As for the interpretation of the control problem, let us just point out 
the following comments:
\begin{itemize}
\item[(i)]
	The cost functional \eqref{costfunct} is designed to track the state variables during 
	the evolution. The targets $\phQ,\sQ,\phO, \sO$, especially $\phO$ and $\sO$, have to be chosen
	as a desirable configuration for clinical reasons, e.g., for surgery.
	Moreover, if some stable configuration for the system has known, it can be taken as well as
	a target.
\item[(ii)] 
	The smaller $\norma{\ph - \phQ}_{L^2(Q)}^2$ is, the closer the solution $\ph $
	is to the target $\phQ$, and the same goes for the other variables.
	On the other hand, the term $\norma{u}_{L^2(Q)}^2$ penalizes the 
	large values of the control variable and it can be 
	read as the side-effect that may occur if too many drugs are dispensed
	to the patient. 
\item[(iii)]
	The ratios between the constants 
	$\bz, \bQ,\bO, \bQh,\bOh$ implicitly describe 
	which targets hold the leading part in the application.
\end{itemize}

Let us anticipate that for our purpose, we have to restrict
the analysis to the case in which the function $P(\ph)$
degenerates to a positive constant $P$. 
Thus, the system we are going to face is the following
\begin{align}
   \a \dt \m_{\a} + \dt \ph_{\a} - \Delta \m_{\a} &=P (\sigma_{\a} - \mu_{\a})
  \quad \hbox{in $\, Q$}
  \label{EQaprima}
  \\
   \mu_{\a} &= \beta \dt \phi_{\a} - \Delta \phi_{\a} + F'(\phi_{\a})
  \label{EQaseconda}
  \quad \hbox{in $\,Q$}
  \\
   \dt \sigma_{\a} - \Delta \sigma_{\a} &= - P (\sigma_{\a} - \mu_{\a}) + u_{\a}
  \label{EQaterza}
  \quad \hbox{in $\,Q$}
  \\
 \dn \m_{\a}  & = \dn \ph_{\a} =\dn \s_{\a} = 0
  \quad \hbox{on $\,\Sigma$}
  \label{BCEQa}
  \\
   \m_{\a}(0)&=\m_0,\, \ph_{\a}(0)=\ph_0,\, \s_{\a}(0)=\s_0
  \quad \hbox{in $\,\Omega,$}
  \label{ICEQa}
\end{align}
\Accorpa\EQa EQaprima ICEQa
where we have written $\m_\a,\ph_\a$ and $\s_\a$ for the state variables to stress
that they are solution to the system in which $\a>0$. 
Such a state system leads to the following control problem:
\Bsist
	\non
	\CP_{\boldsymbol{\a}} && \hbox{Minimize $\J(\ph,\m,u)$}
	\hbox{ subject to the control contraints \eqref{Uad} and under the}
	\\ && \non
	\hbox{requirement that the variables $(\ph, \s)$ 
	solve the system \EQa.}
	\label{CPa}
\Esist
On the other hand,
we will denote with the symbols $\m,\ph$ and $\s$ their corresponding limits as $\a\searrow 0$.
The asymptotic behavior of the above system, as $\a \searrow 0$,
has been one of the main features of \cite{CGRS_ASY}. More precisely, in 
\cite[Thm.~2.5 and Thm.~2.6]{CGRS_ASY} the authors discuss
the passage to the limit as $\a \searrow 0$ and
rigorously proved in which sense system \EQa~converge to 
\begin{align}   
  \dt \ph - \Delta \m 
  &= P (\sigma - \mu)
  \quad \hbox{in $\, Q$}
  \label{EQprima}
  \\
   \mu &= \beta \dt \phi - \Delta \phi + F'(\phi)
  \label{EQseconda}
  \quad \hbox{in $\,Q$}
  \\
   \dt \sigma - \Delta \sigma & = - P(\sigma - \mu) + u
  \label{EQterza}
  \quad \hbox{in $\,Q$}
  \\
   \dn \m &= \dn \ph =\dn \s = 0
  \quad \hbox{on $\,\Sigma$}
  \label{BCEQ}
  \\
   \m(0)&=\m_0,\, \ph(0)=\ph_0,\, \s(0)=\s_0
  \quad \hbox{in $\,\Omega,$}
  \label{ICEQ}
\end{align}
\Accorpa\EQ EQprima ICEQ
having the care of showing the restrictions under 
which existence and uniqueness hold, respectively.
Moreover,
they also exhibit an error estimate between the solution to system 
\EQa~and the solution to \EQ, which in turn implies the uniqueness 
to the second. Note that to address the corresponding 
control problem, the uniqueness of system \EQ~is mandatory.

Therefore, the control problem we want to solve in this paper 
can be summarized as follows:
\Bsist
	\non
	\CP && \hbox{Minimize $\J(\ph,\m,u)$}
	\hbox{ subject to the control contraints \eqref{Uad} and under the}
	\\ && \non
	\hbox{requirement that the variables $(\ph, \s)$ 
	yield a solution to \EQ.}
	\label{CP}
\Esist
Here, let us sketch some
strategies which are usually employed in control theory 
for the class of linear-quadratic 
cost functional, referring to, e.g., \cite{Trol, Lions_OPT} 
for a complete and thorough presentation.
The main aim of control theory is to prove the existence (eventually also 
uniqueness) of optimal control and provide some necessary
(and eventually sufficient) conditions for optimality.
Once that the well-posedness of the state system has been proved,
the existence of optimal controls easily follows by combining the lower weak sequential
semicontinuity of the cost functional $\J$ with standard weak compactness results for
reflexive Banach spaces. On the other hand, in the nonlinear constrained PDEs control theory,
usually, the uniqueness is out of reach. In fact, ordinarily, one appeal to the strict convexity 
of the cost functional to infer uniqueness from the existence part,
but, whenever the state system, and therefore the corresponding control-to-state
operator, is nonlinear, one cannot hope to recover the strict convexity.
The second step consists of looking for
some optimality conditions. As a matter of fact, since the set of admissible controls
is convex, it follows from standard results of convex analysis that the necessary
condition for optimality
of $\overline{u}$ is carried out by the following variational inequality
\Beq
	\label{optimal_formal}
	D \J_{red}(\overline{u})(v-\overline{u}) \geq 0 \quad \hbox{for every $v \in \Uad$,}
\Eeq
where $D\J_{red}$ stands for the derivative of the reduced cost functional in a 
proper functional sense. Moreover, let us recall that 
$\J_{red}$ is essentially obtained as a composition of 
the cost functional $\J$ and the control-to-state operator $\S$. 
So, since $\J$ is trivially \Frechet~differentiable, the classical 
technique relies on proving the \Frechet~differentiability of $\S$ and
then invoke the chain rule to conclude. 
Anyhow, the above procedure does not lead to the desired conclusion
since does not provide an explicit characterization of the gradient
$\nabla J_{red}(\overline{u})$.
Hence, as in the classical constrained control theory, the Lagrange multipliers can be 
introduced to include the constraints in the minimization problem.
This requires to solve another system,
called adjoint problem and which variables are called adjoint, or 
co-state, variables. Finally, after solving this latter, the variational inequality
\eqref{optimal_formal} can be expressed in a more convenient
way which directly allows
us to represent $\nabla J_{red}(\overline{u})$.
To conclude this overview of control theory, let us emphasize that the second-order
derivative can give us meaningful information for the sufficiency. 
However, this is usually less investigated since it introduces some further
technicalities.
Just to give a simple motivation note that, formally, if a map $\S : {\cal U} \to {\cal Y}$, 
then it follows that $D\S : \cal U \to {\cal L} (\cal U, \cal Y)$ and 
$D^2\S : \cal U \to {\cal L} (\cal U, {\cal L} (\cal U, \cal Y))$.
Hence, if one would like to show that $D\S$ is \Frechet\ differentiable by
checking the definition, it has to consider a double increment 
leading to some technical calculations.
On the other hand, these issues can be overcome by employing some
advanced techniques (see, e.g., \cite{CRT, Trol}).

Summing up, in this paper we aim to show that we can
let the parameter $\a$ goes to zero in $(CP)_\a$ to solve $(CP)$. 
We will provide the classical results for the optimal control; namely,
the existence of optimal control and the first-order necessary condition
for optimality. 
This strategy has a huge advantage. Indeed, we will avoid the 
non-trivial discussion of the \Frechet~differentiability
of the control-to-state mapping corresponding to the state system \EQ.
On the other hand, by adopting this approximation scheme, we need to
overcome an approximation issue since it is not trivially ensured that
every optimal control for $(CP)$ can be approximated by sequences of optimal
controls for $(CP)_\a$.

{\bf Plan of the paper}
We conclude the section by sketching an outline
of the paper. In Section~\ref{SEC_GENERAL_ASS_RESULTS}, we will focus
the attention on two aspects; the first one is setting the framework
and the notation, while the second consists in presenting the
obtained results. In Section~\ref{SEC_EXISTENCE_APPROXIMATION}, 
we start with the corresponding
proofs by checking the existence of optimal control and showing an 
approximation result that will be of crucial importance for the asymptotic 
analysis. In Section~\ref{SEC_OPT_COND}, 
we investigate the asymptotic analysis of the adjoint system 
proving its well-posedness in a proper framework. Lastly, we exploit 
the adjoint problem and the approximation result to provide the first-order 
necessary condition for optimality, reading as a variational inequality.

\section{General Assumptions and Results}
\label{SEC_GENERAL_ASS_RESULTS}
\setcounter{equation}{0}
Let us now come to present the mathematical framework and
state the main results.
First of all, we recall that the set $\Omega$ models 
the tissue where the evolution takes place and we assume it to be an open,
bounded and regular domain in $\erre^3$.
Moreover, for an arbitrary Banach space $X$, we convey to use
$\norma{\cdot}_{X}$ to denote its norm,
the standard symbol $X^*$ for its topological dual, and ${}_{X^*}\<\cdot, \cdot >_X$ 
for the corresponding duality product between $X^*$ and $X$.
Likewise, for every $p \in [1,+\infty]$, 
we use the symbol $\norma{\cdot}_{p}$ for the usual norm in $L^p(\Omega)$.
Since in what follows we are going to use several times some particular spaces, 
it turns out to be convenient to set the following conventions
\Beq
	\non
	H:= \Ldue, \quad V:= \Huno, \quad W:=\graffe{v \in \Hdue : \dn v = 0 \hbox{ on } \Gamma},
\Eeq
where $\partial _n$ stands for the outward normal derivative 
of $\Gamma$, and where
these spaces are equipped with their standard norms in order to have Banach spaces.
Let us remark that the canonical injections 
$V \hookrightarrow H \cong H^* \hookrightarrow V^*$ 
are both continuous and dense. 
Therefore, the triplet $(V,H,\Vp)$ forms a Hilbert triplet. 
Indeed, we can identify, in the usual way,
the duality product of $V$ with the inner product of $H$ as follows
\Beq
	\non
	{}_{\Vp}\< u,v >_V = \iO u v \quad \hbox{for every $u \in H$ and $v \in V$.}
\Eeq

Before diving into the setting we are going to use,
let us underline again that our starting point is the distributed 
optimal control investigated in \cite{S}
which considers a quite strong framework in order to handle
the tricky case of the logarithmic potential.
On the other hand, in order to apply the asymptotic strategy mentioned above,
we have to guarantee the well-posedness of system \EQ~which has been treated in \cite{CGRS_ASY} 
(see also \cite{CGRS_VAN}). 
Hence, the simplification introduced in this second work, in order to treat the
asymptotic analysis, cannot be avoided. 
So, all the results proved in \cite{S} hold since the following setting 
perfectly fits the one there considered.

As the assumptions for the above systems and the cost functional are concerned, 
we require that
\begin{align}
	& 
	\a, \b > 0
	\label{ab}
	\\ &
	\bz, \bQ, b_2, b_3, b_4   \, \hbox{ are nonnegative constants, but not all zero}
	\label{constants}
	\\ &
	\phQ, \sQ \in L^2(Q), \phO, \sO \in \Huno, u_*,u^*\in L^\infty(Q) \hbox{ with } u_*\leq u^* \,\aeQ
	\label{targets}
	\\ &
 	P \hbox{ is a positive constant}
	\label{P}
	\\ &
	\ph_0 \in W, \m_0 \in \Huno, \s_0 \in \Huno.
	\label{initial_data}	
\end{align}
As for the control-box, we assume it to be a closed and convex set, and we 
also owe to the following notation
\begin{align}
	& \non
	\UR \subset L^2(Q) \hbox{ be a non-empty and bounded open set such that
	it contains } \Uad
	\\ \non &
	 \hbox{and } \norma u_2 \leq R \hbox{ for all }
	u \in \UR.
\end{align}
Moreover, as for the nonlinear double-well potential $F$, we postulate that
\Beq
	\label{defF}
	F:\erre \to [0, +\infty), \quad  \hbox{with} \quad  
	F:= \hat{B} + \hat{\pi},
\Eeq
where
\begin{align}
	& 
	\hat{B}:\erre \to [0, +\infty) \ \ \hbox{is convex and lower semicontinuous, }	
	\label{Bhat}
	\hbox{with $\hat{B}(0)=0$}
	\\ &
	\hat{\pi} \in C^1(\erre) \,\, \hbox{is nonnegative,} \,\, \pi:=\hat{\pi}' 
	\hbox{ is \Lip~continuous.}
	\label{pihat}
\end{align}
It follows from the above requirements that $B:=\partial\hat{B}$ is a maximal 
and monotone graph $B\subset \erre \times \erre$ 
(see, e.g., \cite[Ex.~2.3.4, p.~25]{BRZ}) and that $D(\hat{B})=\erre$.
Furthermore, from \eqref{pihat}, we also deduce that $\hat{\pi}$ grows at
most quadratically and that its derivative $\pi$ is linearly bounded.
Unfortunately, to manage the optimal control problem introduced above, we 
are forced to restrict the class of admissible potentials by 
requiring some explicit growth assumptions. In fact, we also assume that
\begin{align}
	\label{F}
	F \hbox{ is a $C^3$ function on $\erre$ satisfying }
	|F''(r)| \leq C_1 (1 + |r|^2),
\end{align}
\Accorpa\tutteleipotesi ab F
for a positive constant $C_1$. 
Despite such a strong framework for the potentials, 
the regular potential \eqref{F_reg} complies with the 
above requirements and can be considered.
We also notice that, owing to
the regularity of the initial datum $\ph_0$ and of \eqref{F},
we realize that $F(\ph_0)\in\Lx1$.
Indeed, from \eqref{F}, we realize that 
$F(r) =O(r^4)$ as $|r| \to +\infty $, and, owing to the Sobolev embedding,
that $\ph_0 \in \Lx4$. 

Now, we start by recalling some already known results. 
First of all, we introduce the well-posedness and the asymptotic results,
as $\a \searrow 0$, for system \EQa. In this regard, we refer to
\cite{CGRS_ASY} and to \cite{CGH,CGRS_VAN}.
As a matter of fact, the following existence result still holds in a rather mild
setting for the potential $F$. Namely, the requirement 
on the potentials can be weakened by assuming that 
$\hat{B}$ may attain also the value $+\infty$
and that for a positive constant $C_B$ it holds that
\Beq
	\non
	|B^\circ(r)| \leq C_B (\hat{B}(r)+1) \quad \hbox{for every $r\in\erre$,}
\Eeq
where $B^\circ(r)$ denotes the element of $B(r)$ having minimum modulus
since, without assuming any regularity property, $B$ may be multivalued.
However, we reinforce the setting according to the uniqueness result 
\cite[Thm.~2.6]{CGRS_ASY} unifying the description, by virtue of simplicity,
as for the control problem both the results are necessary.

In view of what already pointed out, it immediately follows from \cite[Thm.~2.5., Thm.~2.6]{CGRS_ASY}
and \cite{S} the following results.

\Bthm
Let \tutteleipotesi~be fulfilled. Then, system \EQa~admits a unique solution 
$(\m_\a,\ph_\a,\s_\a)$ which satisfies 
\begin{align}
	& 
	\label{reg_statea}
	\m_\a,\s_\a,
	\ph_\a \in \H1 H \cap \L\infty V \cap \L2 W \subset \C0 V.
\end{align}
\Ethm

\Bthm
\label{THM_WELLPOSEDNESS_STATE}
Suppose that \tutteleipotesi~are fulfilled. 
Moreover, for given $\a,\b \in (0,1)$ and $u_\a \in \UR$,
let us denote with
$ (\m_\a,\ph_\a,\s_\a)$ the unique solution to system \EQa~enjoying \eqref{reg_statea}.
Then, there exist $\m,\ph,\s$ and a not relabeled subsequence such that,
as $\a\searrow 0$, we have that
\begin{align}
	\label{ma_to_m}
	\m_\a &\to \m\ \ \hbox{weakly in } \L2 V
	\\ 
	\non
	\ph_\a &\to \ph \ \ \hbox{weakly star in } \H1 H \cap \L\infty V \cap \L2 W 
			\\ & \ \ \  \ \hbox{and strongly in } \L2 V \cap \C0 H
	\label{pha_to_ph}
	\\ 
	\non  
	\s_\a &\to \s \ \ \hbox{weakly star in } \H1\Vp \cap \L\infty H \cap \L2 V
			\\ & \ \ \  \ \hbox{and strongly in } \L2  H	
	\label{sa_to_s}
	\\   
	\label{ama_to_zero}
	\a \m_\a &\to 0 \ \ 
	\hbox{strongly in } \H1 \Vp \cap \L2 V.
\end{align}
\Accorpa\Stateconv  ma_to_m ama_to_zero
Furthermore, there exists a positive constant $K_1$, independent of $\a$, such that
\begin{align}
	\non
	& \a^{1/2}\norma{\m_\a}_{\H1 \Vp}
	+ \norma{\m_\a}_{\L2 V}
	+ \norma{\ph_\a}_{ \H1 H \cap \L\infty V \cap \L2 W }
	\\ & \quad
	+ \norma{\s_\a}_{\H1\Vp \cap  \L\infty H \cap \L2 V}
	\leq K_1.
	\label{stima_state}
\end{align}
In addition, the limit triple $(\m,\ph,\s)$ is the unique solution 
to \EQ~and possesses the following regularity
\begin{align}
	& \m \in \L2 V
	\label{reg_m}
	\\
	& \ph \in \H1 H \cap \L\infty V \cap \L2 W
	\label{reg_ph}	
	\\
	& \s \in \H1\Vp  \cap \L\infty H \cap \L2 V.
	\label{reg_s}
\end{align}
\Ethm
Let us also point out that our assumption perfectly fits the framework of 
\cite{S}. Hence, the following existence result holds.
\Blem
\label{LM_existencePCa}
Assume that assumptions \tutteleipotesi~are fulfilled. Then, for every 
$\a\in(0,1)$, the optimal control problem $(CP)_{\a}$ admits, at least, a solution.
\Elem 
Our next goal is to investigate the asymptotic analysis
for the corresponding adjoint system which has been already studied in \cite{S},
and reads as follows
\begin{align}
  & \b \dt q_{\a} - \dt p_{\a} + \Delta q_{\a} - F''(\bph_{\a})q_{\a} 
  = \bQ(\bph_{\a} - \phQ)
  \quad \hbox{in $\, Q$}
  \label{EQAggaprima}
  \\
  & q_{\a} -\a \dt p_{\a} - \Delta p_{\a} + P(p_{\a} -r_{\a}) = 0
  \label{EQAggaseconda}
  \quad \hbox{in $\,Q$}
  \\
  & -\dt r_{\a} - \Delta r_{\a} + P(r_{\a} - p_{\a}) = \bQh(\bs_{\a} - \sQ)
  \label{EQAggaterza}
  \quad \hbox{in $\,Q$}
  \\[0.1cm]
  &\dn q_{\a} = \dn p_{\a} = \dn r_{\a} = 0
  \quad \hbox{on $\Sigma$}
  \label{BCEQAgga}
  \\ 
  & p_{\a}(T) - \b_{\a} q(T) =\bO(\bph_{\a}(T) - \phO),  \quad
  \a p_{\a}(T) = 0, 
  \quad r_{\a}(T) = {\bOh (\bs_{\a}(T) - \sO)}
  \quad \hbox{in $\Omega$.}
  \label{ICEQAgga}
\end{align}
\Accorpa\EQAgga EQAggaprima ICEQAgga
Again, as a consequence of \cite[Thm.~2.8]{S}, we have at disposal the following result.
\Bthm
\label{THM_existence_adjoint_a}
Assume that the assumptions \tutteleipotesi~are fulfilled. Then, there exists a unique triplet 
$(q_\a,p_\a,r_\a)$ which solves \EQAgga~and possesses the beneath regularity
\Bsist
	\label{reg_adj_a}
	q_\a , p_\a, r_\a  \in  \H1 H \cap \L\infty V \cap \L2 W
	\subset \C0 V.
\Esist
\Ethm
Next, starting from this, one can obtain the first-order necessary condition 
for optimality (see \cite[Thm.~2.9]{S}).
\Bthm
\label{THM_foc_finale_a}
Assume that \tutteleipotesi~are verified. Let $\overline{u} \in \Uad$ be an optimal control for $(CP)_{\a}$,
and let $(\bm_{\a},\bph_{\a},\bs_{\a})$ and $(p_{\a},q_{\a},r_{\a})$
be the corresponding optimal state and co-state, respectively. 
Then, it follows that
\Beq
  \label{foc_finale_a}
  \intQ (r_{\a} + \bz \overline{u}_{\a})(v - \overline{u}_{\a})
  \geq 0 \quad \forall v \in \Uad.
\Eeq
\Ethm
%
Now, let us state the novelties.
We aim at showing that, as $\a \searrow 0$, the above system converge, in a proper
sense, to the adjoint system corresponding to \EQ~which reads as
\begin{align}
  & \b \dt q - \dt p + \Delta q - F''(\bph)q = \bQ(\bph - \phQ)
  \quad \hbox{in $\, Q$}
  \label{EQAggprima}
  \\
  & q - \Delta p + P(p -r) = 0
  \label{EQAggseconda}
  \quad \hbox{in $\,Q$}
  \\
  & -\dt r - \Delta r + P(r - p) = {\bQh (\bs - \sQ)}
  \label{EQAggterza}
  \quad \hbox{in $\,Q$}
  \\[0.1cm]
  &\dn q = \dn p = \dn r = 0
  \quad \hbox{on $\Sigma$}
  \label{BCEQAgg}
  \\ 
  & p(T) - \b q(T) =\bO(\bph(T) - \phO),  \quad
  \a p(T) = 0, 
  \quad r(T) = {\bOh (\bs(T) - \sO)}
  \quad \hbox{in $\Omega$.}
  \label{ICEQAgg}
\end{align}
\Accorpa\EQAgg EQAggprima ICEQAgg
We claim that, under suitable assumptions, the above system admits a unique
solution in a variational sense. 
To avoid ambiguity, let us introduce the notion of solution we are going
to employ for this latter.
\Bdef
\label{DEF_solution}
The triplet $(q,p,r)$ is a solution to system \EQAgg~if
it satisfies the variational formulation
\begin{align*}
	& \non
	- {}_{\Vp }\< \dt (p-\b q)(t) , v >_{ V}
	-  \iO \nabla q(t) \cdot \nabla v
	-  \iO F''(\bph(t))q(t) v
	\\ & \hspace{5.2cm}
	= \iO \bQ(\bph(t) - \phQ(t))v
	\quad \hbox{for every $v\in V,\ \aat$}
	\\ \non & \quad
	\iO q(t) v
	+ \iO \nabla p(t) \cdot \nabla v 
	+ P \iO (p(t) - r(t)) v 
	= 0
	\quad \hbox{for every $v\in V,\ \aat$}
  	\\ \non
  	& - {}_{\Vp}\< \dt r(t), v >_V
  	+ \iO \nabla r(t) \cdot \nabla v 
  	+ P \iO (r(t) - p(t)) v
  	\\ & \hspace{5.2cm}
  	=  \iO \bQh(\bs(t) - \sQ(t)) v
  	\quad \hbox{for every $v\in V, \ \aat$,}
\end{align*}
and the final conditions
\Bsist
	\non
	\iO (p-\b q)(T) v 
	= 
	 \iO \bO(\bph(T) - \phO) v
	\quad \hbox{for every $v \in V$}
\Esist
and
\Bsist
	\non
	\iO r(T) v 
	= 
	\iO \bOh (\bs(T) - \sO) v
	\quad \hbox{for every $v \in V$.}
\Esist
Moreover, it has to possess the following regularity
\begin{align}
	\label{regq}
	 q &\in \L2 V
	\\ 
	\label{regp}
	p & \in \L2 W
	\\   
	\label{regr}
	 r &\in \H1\Vp \cap \L\infty H \cap \L2 V \subset \C0 H
	\\   
	\label{regpq}
	 p-\b q &\in \H1\Vp \cap \L\infty H \cap \L2 V \subset \C0 H.
\end{align}
\Accorpa\RegAgg regq regpq
\Edef
Thus, we are in a position to introduce the results concerning the
asymptotic behavior of system \EQAgga, which will be fundamental for the asymptotic 
investigation.
\Bthm
\label{THM_well-posedness_adjointsystem}
Assume that \tutteleipotesi~are in force. 
Let $(q_\a, p_\a, r_\a)$ be the unique solution to \EQAgga~satisfying \eqref{reg_adj_a}. 
Then, as $\a \searrow 0$, and up to a not relabeled subsequence, we have that
\begin{align}
	\label{q_atoq}
	q_\a &\to q\ \ \hbox{weakly in } \L2 V
	\\ 
	\label{p_atop}
	p_\a &\to p \ \ \hbox{weakly in } \L2 W
	\\   
	\label{r_ator}
	r_\a &\to r \ \ \hbox{weakly star in } \H1\Vp \cap \L\infty H \cap \L2 V
			\\ & \ \ \  \ \hbox{and strongly in } \L2 H
	\\
	p_\a - \b q_\a &\to p - \b q \ \ \hbox{weakly star in } \H1\Vp \cap \L\infty H \cap \L2 V
			\\ & \ \ \  \ \hbox{and strongly in } \L2 H
	\\	
	\label{apa_0}
	\a p_\a &\to 0 \ \ \hbox{weakly star in } \H1 H \\ & \ \ \  \ 
		\hbox{and strongly in } \L\infty V\cap \L2 W.
\end{align}
\Accorpa\convadjab q_atoq apa_0
Moreover, there exists a positive constant $K_2$, independent of $\a$, such that
\begin{align}
	& \non
	\norma{p_\a - \b q_\a}_{\H1 \Vp \cap \L\infty H \cap \L2 V}
	+\norma{q_{\a}}_{\L2 V}
	+\a^{1/2} \norma{p_{\a}}_{\L\infty V}	
	\\ & \qquad
	+ \a \norma{p_{\a}}_{\H1 H}	
	+ \norma{p_{\a}}_{\L2 W}
	+ \norma{r_{\a}}_{\H1 \Vp \cap \L\infty H \cap \L2 V}
	\leq K_2.
\end{align}
In addition, the limit $(q,p,r)$ is the unique solution to problem \EQAgg\
in the sense of Definition \ref{DEF_solution}.
\Ethm

With all these results at disposal, we can announce
the results regarding the existence of optimal controls and the first-order 
necessary condition that every optimal control has to satisfy.
\Bthm
\label{THM_existence_of_opt_control}
Suppose that \tutteleipotesi~are satisfied. 
Then, the optimal control problem $(CP)$ admits, at least, a solution
$\overline{u} \in \Uad$.
\Ethm

\Bthm
\label{THM_foc_finale}
Assume that \tutteleipotesi~are in force and let $\overline{u} \in \Uad$ 
be an optimal control for $(CP)$
with its corresponding optimal state $(\bm,\bph,\bs)$. Moreover, let us denote by 
$(p,q,r)$ the associated solution to the adjoint system \EQAgg. 
Then, the necessary condition for optimality of $\overline{u}$ is given by
the following variational inequality 
\Beq
  \label{foc_finale}
  \intQ (r + \bz \overline{u})(v - \overline{u})
  \geq 0 \quad \forall v \in \Uad.
\Eeq
Furthermore, whenever $\bz \not= 0$, the optimal control $\overline{u}$ is the 
$\L2 H-$projection of $-{r}/{\bz}$ onto the closed subspace $\Uad$.
\Ethm
Let us emphasize a consequence which is of \sfw~importance for the numerical approach.
Comparing the expected theoretical condition \eqref{optimal_formal}
with the explicit \eqref{foc_finale},
it immediately follows that we can identify, via Riesz's representation
theorem, the gradient of the reduced cost functional
as $\nabla \J_{red}(\overline{u})=r + \bz \overline{u}$.
Hence, for the numerical approach, the optimal control problem can be viewed
as a constrained minimization of a function, $\J_{red}$, of which we know
the gradient (think of the \wk~projected conjugate gradient method).

In the remainder of the section, we recall some \wk~results which will be useful 
later on. At first, let us remind the Young inequality
\Beq
  ab \leq \delta a^2 + \frac 1 {4\delta} \, b^2
  \quad \hbox{for every $a,b\geq 0$ and $\delta>0$}.
  \non
\Eeq
In addition, we often owe to the standard Sobolev embedding
\Beq
	\label{VinL6}
	\Huno \hookrightarrow L^q(\Omega)
	\quad \hbox{which holds for every $q\in [1,6]$.}
\Eeq

In the whole of the paper, let us convey to use the symbol small-case
$c$ for every constant
which only depend on the structural data of the problem, that is:
on the final time~$T$, on~$\Omega$, on $R$, on the shape of the 
nonlinearities, on the norms of the involved functions, and possibly on $\b$. 
Differently, the capital letters are devoted to indicating some specific constant
which we eventually will refer in the sequel. Moreover, since we aim to
let $\a \searrow 0$, we will keep track at every step of the eventual
dependence of the appearing constants by $\a$.

\section{Existence and Approximation of Optimal Controls}
\label{SEC_EXISTENCE_APPROXIMATION}
\setcounter{equation}{0}
\subsection{Existence of Optimal Controls}

Here, we check the existence of optimal controls by proving 
Theorem~\ref{THM_existence_of_opt_control}.
\proof[Proof of Theorem \ref{THM_existence_of_opt_control}]
The method we are going to employ is the celebrated direct methods of calculus of variations.
To begin with, let us pick an arbitrary sequence $\graffe{\a_n}_n \subset (0,1]$ which goes
to zero as $n\to \infty$. Then, we take as $\graffe{u_n}_n:=\graffe{u_{\a_n}}_n$
a minimizing sequence for the cost functional $\J$ constitutes by elements of $\Uad$,
which, for every $n$, is optimal controls for $(CP)_{\a_n}$, which
exist by virtue of Lemma \ref{LM_existencePCa}.
Next, at every step, 
we introduce $(\m_n,\ph_n,\s_n)$ as the solution associated
to system \EQa~with $u=u_n$.
By recalling estimate \eqref{stima_state} and the boundedness of $\Uad$,
it \sfw ly follows from standard weak and weak-star compactness results that
there exists a subsequence, which we do not relabel, some $\overline{u} \in \Uad$ and a triplet 
$(\bm,\bph,\bs)$ such that, as $n \to \infty$, we have that
\Bsist
	&& \non
	u_n \to \overline{u} \ \ \hbox{weakly star in } L^\infty(Q)
	\\ \non &&
	\m_n \to \bm  \ \ \hbox{weakly in } \L2 V
	\\ \non && 
	\ph_n \to \bph \ \ \hbox{weakly star in } \H1 H \cap \L\infty V \cap \L2 W 
	\\ \non &&
	\s_n \to \bs \ \ \hbox{weakly star in } \H1 \Vp \cap \L\infty H \cap \L2 V.
\Esist
Moreover, compactness arguments (see, e.g., \cite[Sec.~8, Cor.~4]{Simon}) 
also yield that
\Beq
	\non
	\ph_n \to \bph \quad \hbox{strongly in }  \C0 H \cap \L2 V,
\Eeq
which gives meaning to the initial condition $\bph(0)=\ph_0$.
This, along with the growth assumption \eqref{F}, 
allows us to infer that
\begin{align}
	\non
	F'(\ph_n)\to F'(\bph) \quad \hbox{strongly in $\L2 H$.}
\end{align}
Then, it suffices to take into account the variational formulation of 
system \EQa, written for $(\m_n,\ph_n,\s_n)$, and pass to the limit as $n \to \infty$.
From this passage, we infer that $(\bm,\bph,\bs)$
is admissible for $(CP)$,
that is $\bm$ and $\bph$ are the unique solutions to \EQ~associated with $\overline{u}$.
Lastly, the weak sequential lower semicontinuity of the cost functional leads to conclude 
that $\overline{u}$ is a minimizer for $(CP)$. 
\qed

\subsection{Approximation of Optimal Controls}
\label{SUB_APPROX}
After existence has been shown, we would like to infer some information
on the behavior of the optimal controls, pointing out
some necessary conditions for optimality. 
We would like to achieve this goal by letting $\a \searrow 0$
in the necessary condition for $(CP)_\a$
expressed by the variational inequality \eqref{foc_finale_a}. 
Although from a formal perspective it could seem reasonable,
we cannot directly proceed this way. In fact,
if we want to let $\a \searrow 0$ without any
restriction, we have to ensure that every optimal control for $(CP)$
can be approximated by a sequence of optimal controls for $(CP)_\a$.
Unfortunately, we are unable to prove such a strong global approximation result.
Anyhow, a partial one can be stated localizing the problem by following the idea
firstly introduced by Barbu in {\cite{BARBU}}. Let us refer
the interested reader, among others, to the contributions \cite{S_DQ,CFGS,CFS,CGS_DQ}, where 
an application of such a technique can be found.
The key ingredient relies on a local perturbation of the cost functional $\J$.
Then, instead of looking for approximating sequence made up by optimal controls for 
$(CP)_\a$, we seek for a sequence of optimal controls for a modified optimization problem.
Namely, we still consider the same state system, whereas we are going to 
minimize the \socal~adapted cost functional which, for every optimal control 
$\overline{u}$ for $(CP)$, is defined by
\Beq
	\label{adapted_cost}
	\widetilde{\J} (\ph,\s,u):=
	\J (\ph,\s,u)
	+ \frac 12 \norma{u-\overline u}_{L^2(Q)}^2.
\Eeq
Due to the fact that the state system is the same, it is \sfw~to
deduce that this slight modification of $\J$ do not change the corresponding 
adjoint system. 
Hence, it is natural to consider the following new minimization
problem:
\Bsist
	\non
	(\boldsymbol{\widetilde{CP})_{\a}} && \hbox{Minimize $\widetilde{\J}(\ph,\m,u)$}
	\hbox{ subject to the control constraints \eqref{Uad} and under the}
	\\ && \non
	\hbox{requirement that the variables $(\ph, \s)$ 
	yield a solution to \EQa.}
	\label{Ptilde}
\Esist
It is worth emphasizing that $\widetilde{\J}$ reduces to $\J$ whenever 
it is restricted to act 
on optimal controls for $(CP)$.
Moreover, the above control problem perfectly complies with the framework of \cite{S} and therefore,
we also have the following lemma at disposal.
\Blem
\label{LM_existencePCtil}
Assume that assumptions \tutteleipotesi~are satisfied. Then, for every 
$\a\in(0,1)$, there exists at least an optimal control for $(\widetilde{CP})_{\a}$.
\Elem
In a similar fashion as above, it also follows from \cite{S}
how the first-order necessary condition for optimality can be outlined:
\Bthm
\label{THM_secondanecab_adapted}
Assume that \tutteleipotesi~are in force. Let $\overline{u}_\a \in \Uad$ 
be an optimal control for $(\widetilde{CP})_{\a}$, and $(\bm_{\a},\bph_{\a},\bs_{\a})$
and $(p_{\a},q_{\a},r_{\a})$ be the corresponding state and co-state,
respectively. Then, the first-order necessary condition
for optimality is characterized by the variational formulation
\Beq
  \label{foc_finalea_adapt}
  \intQ (r_{\a} + \bz \overline{u}_\a + (\overline u_\a - \overline u ))(v - \overline{u}_\a)
  \geq 0 \quad \forall v \in \Uad.
\Eeq
\Ethm
With all these ingredients, we are finally in a position to introduce the aforementioned 
approximation result.
\Bthm
\label{THM_Approximation}
Assume that \tutteleipotesi~are fulfilled. Let us denote $(\bph,\bs,\overline u)$ 
an optimal triplet for $(CP)$ and
let $\graffe{\a_n}_n\subset (0,1]$ be a sequence which goes to zero as $n\to \infty$.
Then, there exists an approximating optimal sequence, namely a sequence
that, for every $n$, consists of an optimal triplet $(\bph_{\a_n}, \bs_{\a_n}, \overline{u}_{\a_n} )$
for $(\widetilde{CP})_{\a_n}$,
such that the following convergences are satisfied
\begin{align}
	 \label{approx_u}
	\overline u_n:=\overline{u}_{\a_n} &\to \overline u \ \ \hbox{strongly in } L^2(Q)
	\\ 
	 \label{approx_ph}
	\bph_n:=\bph_{\a_n}  &\to \bph \ \ \hbox{weakly star in }  \H1 H \cap \L\infty V \cap \L2 W
	\\   \label{approx_b}
	\bs_n:= \bs_{\a_n} &\to \bs \ \ \hbox{weakly star in }  \H1 \Vp \cap \L\infty H \cap \L2 V
	\\  \label{approx_contr}
	\widetilde{\J} (\bph_n, \bs_n, \overline{u}_n ) & \to \J(\bph,\bs,\overline u)
\end{align}
\Accorpa\approx approx_u approx_contr
as $n \to \infty$, and up to a not relabeled subsequence.
\Ethm
\Bdim
By virtue of Lemma \ref{LM_existencePCtil}, for every $n\in\enne$, we can take
an optimal triplet $(\bph_{\a_n}, \bs_{\a_n}, \overline{u}_{\a_n})$ for $(\widetilde{CP})_{\a_n}$ 
that, for convenience, we will denote by $(\bph_n, \bs_{n}, \overline{u}_{n})$.
From the bound pointed out by estimate \eqref{stima_state}, together with the
boundedness of the control-box, after extraction of a subsequence, we easily get that
\Bsist
	&&
	\non
	\overline{u}_n \to {u} \ \ \hbox{weakly star in } L^\infty(Q)
	\\  &&\non
	\bph_n\to \ph \ \ \hbox{weakly star in }  \H1 H \cap \L\infty V \cap \L2 W
	\\  &&\non
	\bs_n \to \s \ \ \hbox{weakly star in }  \H1 \Vp \cap \L\infty H \cap \L2 V.
\Esist
On the other hand, the continuity of the control-to-state mapping entails that
the limit triplet $(\m,\ph,u)$ is admissible for the control problem $(CP)$, that is
$\ph$ and $\s$ are the solution to \EQ~corresponding to $u$. Thus, our purpose is 
now checking that the limit $u$ is not only admissible, but it is actually 
optimal which, in turn,
will imply that $\ph$ and $\s$ are the corresponding optimal states.
In this direction, we rely on monotonicity arguments.
Firstly, the optimality of $(\bph_{n}, \bs_{n}, \overline{u}_{n})$
for $(\widetilde{CP})_{\a_n}$ yields that
\Beq
	\non
	\widetilde{\J} (\bph_n, \bs_n, \overline u_n )
	\leq \widetilde{\J} (\bph,\bs,\overline u)
	\quad \hbox{ for every $n\in\enne$}
\Eeq
and passing to the superior limit to both sides and exploiting 
the definition of the adapted cost functional $\widetilde{\J}$, we realize that
\Beq
	\label{limsup}
	\limsup_{n\to\infty} \widetilde{\J} (\bph_n, \bs_n, \overline u_n )
	\leq
	\widetilde{\J} (\bph,\bs,\overline u)
	= \J (\bph,\bs,\overline u).
\Eeq
Moreover, it is \sfw~to see that also $\widetilde{\J}$ is 
weak sequential lower semicontinuous,
which implies that
\begin{align}
	\non 
	\liminf_{n\to\infty} \widetilde{\J} (\bph_n, \bs_n, \overline{u}_n )
	& \geq \widetilde{\J} (\ph,\s,u) 
	= \J (\ph,\s,u) +\frac 12 \norma{u-\overline{u}}^2_{L^2(Q)}
	\\ & 	
	\label{liminf}
	\geq \J (\bph,\bs,\overline u) +\frac 12 \norma{u-\overline u}^2_{L^2(Q)},
\end{align}
where the optimality of $(\bph,\bs, \overline u)$ for $(CP)$ and the definition of
the adapted cost functional have been invoked.
By combining \eqref{limsup} with \eqref{liminf}, we get the first convergence
we are looking for since it follows that we have arrived at the identity
\Beq
	\frac 12 \norma{u- \overline u}^2_{L^2(Q)} = 0,
	\label{andrea1}
\Eeq
so that $\overline{u}_n $ weakly star converges to $\uopt$.
These limits, also lead us to infer that the triplet $(\ph,\s,u)$ is nothing but
$(\bph,\bs,\overline{u})$.
As \eqref{approx_contr} is concerned, it suffices to remember that 
$\J(\bph,\bs,\overline u) = \widetilde{\J} (\bph,\bs,\overline u)$ and the fact
that the above inferior and superior limits coincide. In fact, 
we realize that the following chain of equality has been shown:
\begin{align*}
	& 
	\lim_{n\to\infty} \widetilde{\J} (\bph_n, \bs_n, \overline u_n )
	= 	\liminf_{n\to\infty} \widetilde{\J} (\bph_n, \bs_n, \overline{u}_n )
	= \limsup_{n\to\infty} \widetilde{\J} (\bph_n, \bs_n, \overline u_n )
	= \J(\bph,\bs,\overline u).
\end{align*}
Thus, we are reduced to prove \eqref{approx_u}. 
Using the above estimates, we infer that
\begin{align}
	& \label{andrea3}
	\J(\bph,\bs,\overline u) 
	= \lim_{n\to\infty} {\J} (\bph_n, \bs_n, \overline u_n )
	+ \frac 12 \norma{\overline u_n - \uopt}^2_{L^2(Q)}.
\end{align}
On the other hand, the lower semicontinuity of the cost functional, along with
the above estimates, entails that 
\begin{align*}
		 \non
		\J(\bph,\bs,\overline u) 
		&\leq \liminf_{n\to\infty} \J(\bph_n, \bs_n, \overline u_n) 
		\leq \limsup_{n\to\infty} \J(\bph_n, \bs_n, \overline u_n) 
		\\  
		&\leq \limsup_{n\to\infty} \widetilde{\J} (\bph_n, \bs_n, \overline u_n) 
		= \lim \widetilde{\J} (\bph_n, \bs_n, \overline u_n) 
		= \J(\bph,\bs,\overline u),
\end{align*}
so that
\begin{align}
	\non
	\J(\bph,\bs,\overline u)
	=
	\lim_{n\to\infty} \J(\bph_n, \bs_n, \overline u_n)
\end{align}
is verified. Therefore, by combining the above property with \eqref{andrea3}
we deduce that
\begin{align*}
	\frac 12 \norma{\overline{u}_n  - \uopt}^2_{L^2(Q)} \to 0,
\end{align*}
which conclude the proof.
\Edim

\section{Optimality Conditions}
\label{SEC_OPT_COND}
\setcounter{equation}{0}
Next, we establish the necessary condition that an optimal control
has to verify. As explained above, in order to pass to the limit in the variational 
inequality \eqref{foc_finale_a}, we have to deal with
the asymptotic of system \EQAgga~and with the approximation issue presented above.
Since the approximating system has been already investigated in the above section,
only the asymptotic analysis of the adjoint system \EQAgga~has been left unanswered.

\subsection{The Adjoint System}
Below, we proceed formally by only providing some a priori estimates.
The justification can be carried out within a Faedo-Galerkin 
scheme as already made in \cite[Sec.~4.4]{S}. 
Let us just point out that, in the approximation, the duality product is replaced
by the $L^2-$inner product and that the final conditions are replaced
by the corresponding $L^2-$orthogonal projection onto the finite space
spanned by the element of the Galerkin basis.

\proof[Proof of Theorem~\ref{THM_well-posedness_adjointsystem}]
%
%
%
%
%
The estimates we are going to perform in a while are twofold. 
Firstly, within a proper approximation scheme,
they will be the key argument to prove the existence of a solution. 
Secondly, since we will keep track at every step of the dependence 
of the appearing constants by $\a$, they will also be the starting
point to let $\a\searrow 0$ to handle the asymptotic analysis of system \EQAgga.

To begin with, it is convenient to rewrite the system \EQAgga~in a different form.
Let us formally motivate this statement; by considering the vanishing of $\a$,
it is \sfw~to realize that the final condition $\a p_\a = 0$ disappears. 
Moreover, by comparing equation \eqref{EQAggaprima} with the corresponding final condition,
it turns out that the variable $p_\a - \b q_\a$ has to be considered as a
single variable,
since only for such a linear combination the final condition is available.
Moving from this consideration, let us set the following notation
\Beq
	\label{defwa}
	w_{\a} := p_{\a} - \b q_{\a},
\Eeq
which, in turn, implies
\Beq
	\non
	q_{\a} = \frac {p_{\a} - w_{\a}} \b \aand p_{\a} = w_{\a} + \b q_{\a}.
\Eeq
According to the above definitions, we rewrite the above system
in terms of the variable $w_\a$ to obtain the new system
\begin{align}
  &
  - \dt w_{\a} 
  + \frac 1\b \Delta p_{\a}
  -\frac 1\b \Delta w_{\a}
  - \frac 1\b F''(\bph_{\a})p_{\a}
  + \frac 1\b F''(\bph_{\a})  w_{\a} 
  = \bQ(\bph_{\a} - \phQ) \quad \hbox{in $\, Q$} 
  \label{EQanew_prima}
  \\
  & \frac  1 \b p_{\a}
  -\frac 1\b w_{\a}
  -\a \dt p_{\a} - \Delta p_{\a} + P(p_{\a} - r_{\a}) = 0
  \label{EQanew_seconda}
  \quad \hbox{in $\,Q$}
  \\
  & -\dt r_{\a} - \Delta r_{\a} + P(r_{\a} - p_{\a}) = \bQh (\bs_{\a} - \sQ)
  \label{EQanew_terza}
  \quad \hbox{in $\,Q$}
  \\
  &\dn w_{\a}= \dn p_{\a} = \dn r_{\a} = 0
  \quad \hbox{on $\Sigma$}
  \label{EQBCanew}
  \\ 
  & w_{\a}(T) =\bO (\bph_{\a}(T) - \phO),  \ \
  \a p_{\a}(T) = 0,  \  \
  r_{\a}(T) = \bOh (\bs_{\a}(T) - \sO)
  \quad \hbox{in $\Omega$.}
  \label{EQICanew}
\end{align}
\Accorpa\EQanew EQanew_prima EQICanew
Now, we start presenting the estimates.

{\bf First estimate}
In the first place, we multiply equation \eqref{EQanew_prima} by $w_{\a}$, 
\eqref{EQanew_seconda} by $p_{\a} - \Delta p_{\a}$, \eqref{EQanew_terza} by $r_{\a}$
and integrate over $Q_t^T$ and by parts to obtain, upon rearranging the terms, that
\Bsist
	\non &&
	\frac 12 \IO2 {w_{\a}}
	+ \frac 1\b \intQtT |\nabla w_{\a}|^2
	+ \frac {\a }2 \iO ( |p_{\a}(t)|^2 +  |\nabla p_{\a}(t)|^2)
	+ \biggl( \frac 1\b +  P \biggr) \intQtT |p_{\a}|^2
	\\ \non && \qquad
	+ \biggl(  \frac 1\b + P + 1 \biggr) \intQtT |\nabla p_{\a}|^2
	+  \intQtT |\Delta p_{\a}|^2
	+ \frac 12 \IO2 {r_{\a}}
	+ \intQtT |\nabla r_{\a}|^2
	+ P \intQtT |r_{\a}|^2
	\\ \non && \quad
	=
	\frac 12 \iO |\bO (\bph_{\a}(T) - \phO)|^2
	+ \frac 12 \iO |\bOh (\bs_{\a}(T) - \sO)|^2
	+ \intQtT \bQ(\bph_{\a} - \phQ) w_{\a}
	\\ \non && \qquad
	+ \intQtT \bQh (\bs_{\a} - \sQ) r_{\a}
	+ \frac 1\b \intQtT F''(\bph_{\a}) p_{\a} w_{\a}
	- \frac 1\b \intQtT F''(\bph_{\a}) w_{\a}^2
	- \frac 2{\b} \intQtT \Delta p_{\a} w_{\a} 
	\\ \non && \qquad
	+  P  \intQtT  r_{\a} (p_{\a} - \Delta p_{\a})
	+ \frac 1 \b \intQtT w_{\a} p_{\a} 
	+ P \intQtT p_{\a} r_{\a},
\Esist
where we denote the terms on the \rhs~by $I_1,...,I_{10}$, in this order.
Moreover, the integrals on the \lhs~are nonnegative, whereas the ones on the \rhs~can
be bounded as follows.
Using the final conditions \eqref{EQICanew}, assumptions \eqref{constants},
\eqref{targets}, and \eqref{reg_statea},
we infer by the Young inequality that
\Beq
	\non
	|I_1| + |I_2| + |I_3| + |I_4| 
	\leq
	c \intQtT (|w_{\a}|^2 + |r_{\a}|^2) + c.
\Eeq
As for $I_5$ and $I_6$, we recall the growth assumption \eqref{F} and the fact that
$\bph_{\a}$, as a solution to \EQa, verifies estimate \eqref{stima_state}. 
Thus, along with
the \Holder~and Young inequalities, and the standard embedding \eqref{VinL6}, we get that
\begin{align}
	\non 
	|I_5| + |I_6| 
	&\leq
	\frac {C_{1}}\b \intQtT (1 + |\bph_{\a}^2|) p_{\a} w_{\a} 
	+ \frac {C_{1}}\b \intQtT (1 + |\bph_{\a}^2|) w_{\a}^2
	\\ \non & 
	\leq 
	c \inttt (1 + \norma{\bph_{\a}^2}_{3}) \norma{p_{\a}}_{6} \norma{w_{\a}}_{2}
	+ c \inttt (1 + \norma{\bph_{\a}^2}_{3}) \norma{w_{\a}}_{6} \norma{w_{\a}}_{2}
	\\ \non & 
	\leq 
	c \inttt (1 + \norma{\bph_{\a}}^2_{6}) \norma{p_{\a}}_{6} \norma{w_{\a}}_{2}
	+ c \inttt (1 + \norma{\bph_{\a}}^2_{6}) \norma{w_{\a}}_{6} \norma{w_{\a}}_{2}
	\\ \non & 
	\leq 
	c \inttt (1 + \norma{\bph_{\a}}^2_{V}) \norma{p_{\a}}_{V} \norma{w_{\a}}_{H}
	+ c \inttt (1 + \norma{\bph_{\a}}^2_{V}) \norma{w_{\a}}_{V} \norma{w_{\a}}_{H}
	\\ \non & 
	\leq 
	\d \inttt \norma{p_{\a}}^2_{V}
	+ \d  \intQtT |\nabla w_{\a}|^2
	+ \cd \intQtT |w_{\a}|^2,
\end{align}
for a positive constant $\d$ yet to be determined.
Next, invoking once more the Young inequality, we argue that
\Bsist
	\non
	|I_7| 
	\leq
	\frac 1 4 \intQtT |\Delta p_{\a}|^2 
	+ \frac 2 {\b} \intQtT |w_{\a}|^2,
\Esist
and also that
\Bsist
	\non
	|I_8| + |I_9| + |I_{10}| 
	\leq
	3 \d \intQtT |p_{\a}|^2
	+ \frac 14 \intQtT |\Delta p_{\a}|^2
	+ \cd \intQtT (|w_{\a}|^2 + |r_{\a}|^2).
\Esist
Hence, upon collecting all these terms, we realize that it suffices to fix $\d$ small enough.
Namely, we pick $\d=\hat{\d}$ such that 
\Beq
	\non
	\hat{\d} < \min 
	\biggl\{ 
	\frac 1 { \b},
	\frac 14 \biggl(\frac 1 \b+ P  \biggr),
	\frac 1 { \b} + P + 1
	\biggr\}.
\Eeq
Therefore, a Gronwall argument yields that
\begin{align}
	\non
	\norma{w_{\a}}_{\L\infty H \cap \L2 V}
	+\a^{1/2} \norma{p_{\a}}_{\L\infty V}	
	+ \norma{p_{\a}}_{\L2 W}
	+ \norma{r_{\a}}_{\L\infty H \cap \L2 V}
	\leq c,
\end{align}
for a suitable positive constant $c$ independent of $\a$.
Moreover, let us note that
\Beq
	\non
	\norma{\a p_\a}_{\L\infty V}
	\leq
	c \, \a^{1/2}.
	\label{andrea5}
\Eeq

{\bf Second estimate}
Multiplying \eqref{EQanew_prima} by an arbitrary $v\in \L2 V$, integrating
over $Q$ and by parts, and making use of the above bounds, we infer that
\begin{align}
	\non 
	\Big | \intQ \dt w_{\a} \, v \Big |
	&\leq 
	c \norma{\nabla p_{\a}}_{\L2 H} \norma{\nabla v}_{\L2 H}
    + c \norma{\nabla w_{\a}}_{\L2 H} \norma{\nabla v}_{\L2 H}
    \\ \non & \quad
    + c \norma{p_{\a}}_{\L2 H} \norma{v}_{\L2 V}
    + c \norma{w_{\a}}_{\L2 H} \norma{v}_{\L2 V}
	+ c \norma{v}_{\L2 H}
    \\ \non & 
	\leq 
	c \norma{v}_{\L2 V}.
\end{align}
Then, dividing both sides by $\norma{v}_{\L2 V}$ and passing to the superior
limit leads to conclude that
\Beq
	\non
	\norma{\dt w_{\a}}_{\L2 \Vp}
	\leq
	c.
\Eeq

{\bf Third estimate}
By the same token, we employ the above estimates to obtain that
\Beq
	\non
	\norma{\dt r_{\a}}_{\L2 \Vp}
	\leq
	c.
\Eeq

{\bf Fourth estimate}
Lastly, comparison in equation \eqref{EQanew_seconda}, along
with the above estimates, produces
\Beq
	\non
	\norma{\a \dt p_{\a}}_{\L2 H} 
	\leq
	c.
\Eeq
It is now a standard matter to show that 
the above estimates will be sufficient, withing a Galerkin
scheme, to provide the existence of a solution to \EQAgg~which also satisfies \RegAgg.
Furthermore, the existence, together with the linearity of the system, 
also implies its uniqueness.

{\bf Passage to the limit}
Here, we draw some consequences from the aforementioned estimates
checking that, in a proper sense, system \EQAgga~converges to \EQAgg.
Owing to standard weak compactness arguments
it turns out that, up to a not relabeled subsequence, the following convergences hold
\begin{align}
	\non
	w_\a &\to w \ \ \hbox{weakly star in } \H1\Vp \cap \L\infty H \cap \L2 V
	\\ 
	\non
	p_\a &\to p \ \ \hbox{weakly in } \L2 W
	\\   
	\non
	r_\a &\to r \ \ \hbox{weakly star in } \H1\Vp \cap \L\infty H \cap \L2 V.
\end{align}
Moreover, the compact embedding of $\H1 \Vp \cap \L2 V$ into $\C0 H$ 
guarantees that the final data are meaningful and that
\begin{align}
	\non
	w_\a &\to w \ \ \hbox{strongly in } \L2 H
	\\ 
	\non
	r_\a &\to r \ \ \hbox{strongly in } \L2 H
	\\   
	\non	
	\a p_\a &\to 0 \ \ \hbox{weakly star in } \H1 H  \ \hbox{and strongly in } \L\infty V \cap \L2 W.
\end{align}
Hence, by combining the above first and second convergences with the definition
of the auxiliary variable $w_{\a}$ given by \eqref{defwa}, we realize that
\Beq
	\label{qa_q}
	q_\a \to q \ \ \hbox{weakly in } \L2 V.
\Eeq
Therefore, the above convergences implies
that the weak limit of $w_{\a}$ can be identified with
$w = p - \b q$. So, in what follows, we are legitimate to conveniently
interchange the variables $w$ and $q$ as convenience.

Now, let us take into account the variational formulation of system \EQAgga. 
It consists of seeking for a triplet $(w_\a,p_\a,r_\a)$ such that
satisfies the following problem
\begin{align}
	& \non
	- {}_{\Vp }\< \dt w_{\a}(t) , v >_{ V}
	-  \iO \nabla q_{\a}(t) \cdot \nabla v
	-  \iO F''(\bph_{\a}(t))q_{\a}(t) v
	\\ & \non	\hspace{4.2cm}
	= \iO \bQ(\bph_{\a}(t) - \phQ(t))v
	\quad \hbox{for every $v\in V,\ \aat$}
	\\ \non & \quad
	\iO q_{\a} v 
	-\a \iO \dt p_{\a}(t) v 
	+ \iO \nabla p_{\a}(t) \cdot \nabla v 
	+ P \iO ( p_{\a}(t) - r_{\a}(t)) v 
	\\ & \non	\hspace{7.9cm}
	= 0
	\quad \hbox{for every $v\in V,\ \aat$}
  	\\ \non
  	& 
  	- {}_{\Vp }\< \dt r_{\a}(t) , v >_{ V}
  	+ \iO \nabla r_{\a}(t) \cdot \nabla v 
  	+ P \iO (r_{\a}(t) - p_{\a}(t)) v
  	\\ & \non	\hspace{4.3cm}
	=  \iO \bQh(\bs_{\a}(t) - \sQ(t)) v
	\quad \hbox{for every $v\in V,\ \aat$.}
\end{align}
Moreover, owing to the final conditions \eqref{EQICanew},
it has also to verify the final conditions
\Bsist
	\non
	\iO w_{\a}(T) v 
	= 
	 \iO \bO(\bph_{\a}(T) - \phO) v
	\quad \hbox{for every $v \in V$}
\Esist
and
\Bsist
	\non
	 \iO r_{\a}(T) v 
	= 
	\iO \bOh (\bs_{\a}(T) - \sO) v
	\quad \hbox{for every $v \in V$.}
\Esist
By virtue of the above discussion, along with the convergences \Stateconv, 
we would conclude that, as $\a \searrow 0$, the above system converges 
to the following problem:
\begin{align}
	& \non
	- {}_{\Vp }\< \dt (p-\b q)(t) , v >_{ V}
	-  \iO \nabla q(t) \cdot \nabla v
	-  \iO F''(\bph(t))q(t) v
	\\ & \non \hspace{4.2cm}
	= \iO \bQ(\bph(t) - \phQ(t))v
	\quad \hbox{for every $v\in V,\ \aat$}
	\\ \non & \quad
	\iO q(t) v
	+ \iO \nabla p(t) \cdot \nabla v 
	+ P \iO (p(t) - r(t)) v 
	\\ & \non \hspace{7.6cm}
	= 0
	\quad \hbox{for every $v\in V,\ \aat$}
  	\\ \non
  	& 
  	-{}_{\Vp}\< \dt r(t), v >_V
  	+ \iO \nabla r(t) \cdot \nabla v 
  	+ P \iO (r(t) - p(t)) v
  	\\ & \non \hspace{4.2cm}
	=  \iO \bQh(\bs(t) - \sQ(t)) v
	\quad \hbox{for every $v\in V,\ \aat$}
\end{align}
with the corresponding final conditions
\Bsist
	\non
	\iO (p-\b q)(T) v 
	= 
	 \iO \bO(\bph(T) - \phO) v
	\quad \hbox{for every $v \in V$}
\Esist
and
\Bsist
	\non
	\iO r(T) v 
	= 
	 \iO \bOh (\bs(T) - \sO) v
	\quad \hbox{for every $v \in V$.}
\Esist
\dafare{
To do that, we multiply the first system by a regular function 
$\d \in C^\infty_c(0,T)$, integrate over $(0,T)$, and then pass to the limit 
accounting for the above estimates. Thus, since the obtained limit system holds for every
$\d \in C^\infty_c(0,T)$, one finally recover the last system.
}
Anyhow, to prove such a passage, we need to handle the asymptotics of the
nonlinear term $F''(\bph_{\a}) q_\a $.
We claim that, by combining the growth assumptions on the potential 
\eqref{F} with the strong convergence \eqref{pha_to_ph}, it follows that
\Beq
	\label{F''strong}
	F''(\bph_{\a}) \to F''(\bph) \quad \hbox{strongly in } \L2 H.
\Eeq
Therefore, by combining \eqref{qa_q} with \eqref{F''strong},
the nonlinear term can be handled since we have 
\Beq
	\non
	F''(\bph_{\a}) q_\a \to F''(\bph) q \quad \hbox{weakly in $\L2 H$},
\Eeq
and this conclude the proof.
\qed

\subsection{First-order Necessary Condition}
In this last section, we are going to prove Theorem \ref{THM_foc_finale}
which gives us the first-order necessary condition for optimality.
\proof[Proof of Theorem \ref{THM_foc_finale}]
As already mention, we try to recover the first-order necessary condition for the control
problem $(CP)$ via asymptotic techniques by
letting $\a \searrow 0$, in a suitable sense, in the variational 
inequality \eqref{foc_finale_a}. 
The main issue has been already introduced above
and consists in the fact that we have to guarantee that every optimal control for $(CP)$
can be found as a limit of a sequence made up by optimal controls for $(CP)_\a$.
This can be overcome by invoking the investigated approximation
result. In fact, we consider a sequence $\graffe{\a_n}\subset (0,1]$
which goes to zero as $n\to \infty$, and introduce the sequence
$\overline{u}_n:= \overline{u}_{\a_n}$ of optimal controls for $(\widetilde{CP})_{\a_n}$
introduced in Theorem \ref{THM_Approximation}. 
After further extraction of a subsequence $\graffe{\a_{n_k}}$, 
the convergence pointed out by \convadjab~and \approx~allow
us to pass to the limit, as $k\to \infty$,
in \eqref{foc_finalea_adapt} to achieve the necessary condition 
we are looking for.

Finally, the last sentence follows from an application
of the \wk~Hilbert projection theorem, since $\Uad$ is a non-empty, closed and convex subset 
of $\L2 H$.
\qed

Finally, due to the structure of the control-box $\Uad$, in the case of $\bz > 0$,
we can provide an equivalent implicit characterization 
of the optimal control (see, e.g., \cite{Trol}).
\Bcor
Suppose that \tutteleipotesi~and that $\bz>0$. 
Then, the optimal control $\overline{u}$ for $(CP)$ satisfies
\Beq
	\non
	\overline{u}(x,t)=\max \bigl\{ 
	u_*(x,t), \min\graffe{u^*(x,t),-\frac 1{\bz} r(x,t)} 
	\bigr\} \quad \aaQ .
\Eeq
\Ecor

\Brem
From a little investigation of Theorems \ref{THM_WELLPOSEDNESS_STATE} and 
\ref{THM_well-posedness_adjointsystem}, one realize that the requirements
$\phO, \sO \in \Hx1$ and $\ph_0\in W,$ and $ \m_0,\s_0 \in V$ turn out to
be superabundant. In fact, 
for the limit optimal control problem $(CP)$, to be meaningful, it suffices that
\Beq
	\label{dati_weak}
	\phO, \sO \in H \ \ \hbox{and} \ \ \ph_0\in V, \ \m_0,\s_0 \in H.
\Eeq
The framework we have introduced was motivated by the fact that, 
in order to manage $(CP)_\a$, we directly rely on the results of \cite{S}.
Thus, it has been chosen by comparing the framework of
\cite{S} with the additional assumptions introduced in \cite{CGRS_ASY}
to deal with the asymptotic analysis of system \EQa.
Indeed, whenever $\a>0$, both the requirements \eqref{targets} 
and \eqref{initial_data} have to be fulfilled.

In such a perspective, one may wonder if 
the given requirements can be somehow weakened.
One possible way to proceed can be to assume \eqref{dati_weak}
and define some regularizing sequences
\begin{align}
	& 
	\non
	\graffe{{\phO^\a}}_{\a}, \graffe{{\sO^\a}}_{\a} \in \Hx1 
	\\
	& \non
	\graffe{\ph_0^\a}_\a\in W, \ \graffe{\m_0^\a}_\a \in V, \ \graffe{\s_0^\a}_\a \in V
\end{align}
which satisfy, as $\a \searrow 0$, the following strong convergences
\begin{align}
	& 
	\non
	\phO^\a \to \phO, \ \sO^\a \to \sO \ \hbox{ strongly in $H$}
	\\
	& \non
	\ph_0^\a \to \ph_0 \ \hbox{ strongly in $V$,}	
 	\ \hbox{and} \
	\m_0^\a \to \m_0, \s_0^\a \to \s_0 \ \hbox{ strongly in $H$}.
\end{align}
Then, for every $\a \in (0,1)$, the initial conditions in the state system \eqref{ICEQa} 
has to be replaced with the approximated version
\Beq
	\non
	\m_{\a}(0)=\m_0^\a,\, \ph_{\a}(0)=\ph_0^\a,\, \s_{\a}(0)=\s_0^\a  \quad \hbox{in $\,\Omega.$}
\Eeq
Moreover, the cost functional $\J$ has to be substituted by
\begin{align}	
	\non
	\J^\a (\ph, \s, u)   &:= 
	\frac \bQ 2 \norma{\ph - \phQ}_{L^2(Q)}^2
	+\frac \bO 2 \norma{\ph(T)-\phO^\a}_{L^2(\Omega)}^2
	+ \frac \bQh 2 \norma{\s - \sQ}_{L^2(Q)}^2
	\\ & \qquad
	\non
	+ \frac \bOh 2 \norma{\s(T)-\sO^\a}_{L^2(\Omega)}^2
    + \frac \bz 2 \norma{u}_{L^2(Q)}^2,
\end{align}
and the adapted cost with ${\widetilde{\J}}^\a$, defined according to 
\eqref{adapted_cost}. Thus, the new $(CP)_\a$ consists of minimizing
the cost functional $\J^\a$ subject to the control-box constraints $\Uad$,
and under the assumption that $\ph, \s$ are solution to 
this new approximated state system, namely system \accorpa{EQaprima}{BCEQa} coupled
with the above initial data.
It is immediately clear that the corresponding investigation
will became more technical and it is not clear if such an effort is worth
to be pursued.
\Erem

\subsection*{Acknowledgments}
The author would like to thank especially one of the referees 
for the careful reading and the precious suggestions which 
have improved the manuscript.

\vspace{3truemm}

\Begin{thebibliography}{10}
\footnotesize

\bibitem{Agosti}
A. Agosti, P.F. Antonietti, P. Ciarletta, M. Grasselli and M. Verani,
A Cahn-Hilliard-type equation with application to tumor growth dynamics, 
{\it Math. Methods Appl. Sci.}, {\bf 40} (2017), 7598-–7626.

\bibitem{BARBU}
V. Barbu,
Necessary conditions for nonconvex distributed control problems governed by 
elliptic variational inequalities, 
{\it J. Math. Anal. Appl.} {\bf 80} (1981), 566-597.


\bibitem{BRZ}
H. Brezis,
``Op\'erateurs maximaux monotones et semi-groupes de contractions dans les
espaces de Hilbert'', North-Holland Math. Stud. {\bf 5}, North-Holland, Amsterdam, 1973.

\bibitem{CRT}
E. Casas, J.C. de los Reyes and F. Tr\"oltzsch,
Sufficient second-order optimality conditions for semilinear control problems
with pointwise state constraints,
{\it SIAM J. Optim.} {\bf 19(2)}, (2008), 616-643.

\bibitem{CRW}
C. Cavaterra, E. Rocca and H. Wu,
Long-time Dynamics and Optimal Control of a Diffuse Interface Model for Tumor Growth,
{\it Appl. Math. Optim.} (2019), https://doi.org/10.1007/s00245-019-09562-5.

\bibitem{CFGS}
P. Colli, M.H. Farshbaf-Shaker, G. Gilardi and J. Sprekels,
Optimal boundary control of a viscous Cahn–-Hilliard system
with dynamic boundary condition and double obstacle potentials,
{\it SIAM J. Control Optim.} {\bf 53} (2015), 2696-2721.
\bibitem{CFS}
P. Colli, M.H. Farshbaf-Shaker and J. Sprekels,
A deep quench approach to the optimal control of an Allen–-Cahn equation
with dynamic boundary conditions and double obstacles,
{\it Appl. Math. Optim.} {\bf 71} (2015), 1-24.

\bibitem{CGH}
P. Colli, G. Gilardi and D. Hilhorst,
On a Cahn-Hilliard type phase field system related to tumor growth,
{\it Discrete Contin. Dyn. Syst.} {\bf 35} (2015), 2423-2442.

\bibitem{SM}
P. Colli, G. Gilardi, G. Marinoschi and E. Rocca,
Sliding mode control for a phase field system related to tumor growth,
{\it Appl. Math. Optim.} to appear (2019), doi:10.1007/s00245-017-9451-z.
%

\bibitem{CGS_DQ}
P. Colli, G. Gilardi and J. Sprekels,
Optimal velocity control of a convective Cahn-–Hilliard system with double obstacles
and dynamic boundary conditions: a `deep quench' approach.
{\it J. Convex Anal.}, to appear (2018).

\bibitem{CGRS_OPT}
P. Colli, G. Gilardi, E. Rocca and J. Sprekels,
Optimal distributed control of a diffuse interface model of tumor growth,
{\it Nonlinearity} {\bf 30} (2017), 2518-2546.

\bibitem{CGRS_ASY}
P. Colli, G. Gilardi, E. Rocca and J. Sprekels,
Asymptotic analyses and error estimates for a \CH~type phase field system modeling tumor growth,
{\it Discrete Contin. Dyn. Syst. Ser. S} {\bf 10} (2017), 37-54.

%
%
%
%

\bibitem{CGRS_VAN}
P. Colli, G. Gilardi, E. Rocca and J. Sprekels,
Vanishing viscosities and error estimate for a Cahn–Hilliard type phase field system related to tumor growth,
{\it Nonlinear Anal. Real World Appl.} {\bf 26} (2015), 93-108.
%



\bibitem{CLLW}
V. Cristini, X. Li, J.S. Lowengrub, S.M. Wise,
Nonlinear simulations of solid tumor growth using a mixture model: invasion and branching.
{\it J. Math. Biol.} {\bf 58} (2009), 723–-763.

\bibitem{CL}
V. Cristini, J. Lowengrub,
Multiscale Modeling of Cancer: An Integrated Experimental and Mathematical
{\it Modeling Approach. Cambridge University Press}, Leiden (2010).

\bibitem{DFRGM}
M. Dai, E. Feireisl, E. Rocca, G. Schimperna, M. Schonbek,
Analysis of a diffuse interface model of multispecies tumor growth,
{\it Nonlinearity\/} {\bf  30} (2017), 1639.

\bibitem{EK_ADV}
M. Ebenbeck and P. Knopf,
Optimal control theory and advanced optimality conditions for a diffuse interface model of tumor growth 
{\it preprint arXiv:1903.00333 [math.OC],} (2019), 1--34.

\bibitem{EK}
M. Ebenbeck and P. Knopf,
Optimal medication for tumors modeled by a Cahn-Hilliard-Brinkman equation,
{\it preprint arXiv:1811.07783 [math.AP],} (2018), 1-26.

\bibitem{EGAR}
M. Ebenbeck and H. Garcke,
Analysis of a Cahn–-Hilliard–-Brinkman model for tumour growth with chemotaxis.
{\it J. Differential Equations,} (2018) https://doi.org/10.1016/j.jde.2018.10.045.

\bibitem{FGR}
S. Frigeri, M. Grasselli, E. Rocca,
On a diffuse interface model of tumor growth,
{\it  European J. Appl. Math.\/} {\bf 26 } (2015), 215-243. 

\bibitem{FLRS}
S. Frigeri, K.F. Lam, E. Rocca, G. Schimperna,
On a multi-species Cahn-Hilliard-Darcy tumor growth model with singular potentials,
{\it Comm. in Math. Sci.} {\bf  (16)(3)} (2018), 821-856. 

\bibitem{FRL}
S. Frigeri, K.F. Lam and E. Rocca,
On a diffuse interface model for tumour growth with non-local interactions and degenerate 
mobilities,
In {\sl  Solvability, Regularity, and Optimal Control of Boundary Value Problems for PDEs},
P. Colli, A. Favini, E. Rocca, G. Schimperna, J. Sprekels (ed.),
{\it Springer INdAM Series,} {\bf 22}, Springer, Cham, 2017.

%
\bibitem{GARL_1}
H. Garcke and K. F. Lam,
Well-posedness of a Cahn-–Hilliard–-Darcy system modelling tumour
growth with chemotaxis and active transport,
{\it European. J. Appl. Math.} {\bf 28 (2)} (2017), 284-316.
\bibitem{GARL_2}
H. Garcke and K. F. Lam,
Analysis of a Cahn-Hilliard system with non-zero Dirichlet 
conditions modeling tumor growth with chemotaxis,
{\it Discrete Contin. Dyn. Syst.} {\bf 37 (8)} (2017), 4277-4308.
\bibitem{GARL_3}
H. Garcke and K. F. Lam,
Global weak solutions and asymptotic limits of a Cahn-Hilliard-Darcy system modelling tumour growth,
{\it AIMS Mathematics} {\bf 1 (3)} (2016), 318-360.
\bibitem{GARL_4}
H. Garcke and K. F. Lam,
On a Cahn-Hilliard-Darcy system for tumour growth with solution dependent source terms, 
in {\sl Trends on Applications of Mathematics to Mechanics}, 
E.~Rocca, U.~Stefanelli, L.~Truskinovski, A.~Visintin~(ed.), 
{\it Springer INdAM Series} {\bf 27}, Springer, Cham, 2018, 243-264.
\bibitem{GAR}
H. Garcke, K. F. Lam, R. N\"urnberg and E. Sitka,
A multiphase Cahn-Hilliard-Darcy model for tumour growth with necrosis,
{\it Mathematical Models and Methods in Applied Sciences} {\bf 28 (3)} (2018), 525-577.

\bibitem{GARLR}
H. Garcke, K. F. Lam and E. Rocca,
Optimal control of treatment time in a diffuse interface model of tumor growth,
{\it Appl. Math. Optim.} {\bf 78}(3) (2018), {495--544}.
\bibitem{GLSS}
H. Garcke, K.F. Lam, E. Sitka, V. Styles,
A Cahn–-Hilliard-–Darcy model for tumour growth with chemotaxis and active transport.
{\it Math. Models Methods Appl. Sci. } {\bf 26(6)} (2016), 1095-–1148.


\bibitem{OHP}
A. Hawkins, J.T Oden, S. Prudhomme,
General diffuse-interface theories and an approach to predictive tumor growth modeling. 
{\it Math. Models Methods Appl. Sci.} {\bf 58} (2010), 723–-763. 

\bibitem{HDPZO}
A. Hawkins-Daarud, S. Prudhomme, K.G. van der Zee, J.T. Oden,
Bayesian calibration, validation, and uncertainty quantification of diffuse 
interface models of tumor growth. 
{\it J. Math. Biol.} {\bf 67} (2013), 1457–-1485. 

\bibitem{HDZO}
A. Hawkins-Daruud, K. G. van der Zee and J. T. Oden, Numerical simulation of
a thermodynamically consistent four-species tumor growth model, Int. J. Numer.
{\it Math. Biomed. Engng.} {\bf 28} (2011), 3–-24.

\bibitem{HKNZ}
D. Hilhorst, J. Kampmann, T. N. Nguyen and K. G. van der Zee, Formal asymptotic
limit of a diffuse-interface tumor-growth model, 
{\it Math. Models Methods Appl. Sci.} {\bf 25} (2015), 1011-–1043.

\bibitem{Kur}
S. Kurima,
Asymptotic analysis for Cahn-Hilliard type phase field systems related to tumor 
growth in general domains,
{\it Math. Methods in the Appl. Sci} {} (2019), 
https://doi.org/10.1002/mma.5520.


\bibitem{Lions_OPT}
J. L. Lions,
Contr\^ole optimal de syst\`emes gouverne\'s par des equations aux d\'eriv\'ees partielles,
Dunod, Paris, 1968.



\bibitem{Mir_CH}
A. Miranville,
The Cahn-Hilliard equation and some of its variants, 
{\it AIMS Mathematics,} {\bf 2} (2017), 479–-544.

\bibitem{MRS}
A. Miranville, E. Rocca, and G. Schimperna,
On the long time behavior of a tumor growth model,
{\it Commun. Pure Appl. Anal.\/} {\bf 8} (2009) 881-912.
%

\bibitem{MirZel}
A. Miranville, S. Zelik, Attractors for dissipative partial differential equations in bounded and
unbounded domains, in “Handbook of Differential Equations: Evolutionary Equations, Vol. IV”
(eds. C.M. Dafermos and M. Pokorny), Elsevier/North-Holland, 103–200, 2008.

\bibitem{S_b}
A. Signori,
Optimal treatment for a phase field system of Cahn-Hilliard 
type modeling tumor growth by asymptotic scheme. 
{\it Preprint: arXiv:1902.01079 [math.AP]} (2019), 1-28.

\bibitem{S_DQ}
A. Signori,
Optimality conditions for an extended tumor growth model with 
double obstacle potential via deep quench approach. 
{\it Preprint: arXiv:1811.08626 [math.AP]} (2018), 1-25.

\bibitem{S}
A. Signori,
Optimal distributed control of an extended model of tumor 
growth with logarithmic potential. 
{\it Appl. Math. Optim.} (2018), https://doi.org/10.1007/s00245-018-9538-1.

\bibitem{Simon}
J. Simon,
{Compact sets in the space $L^p(0,T; B)$},
{\it Ann. Mat. Pura Appl.\/} 
{\bf 146~(4)} (1987) 65-96.

\bibitem{SW}
J. Sprekels and H. Wu,
Optimal Distributed Control of a Cahn-–Hilliard–-Darcy System with Mass Sources,
{\it Appl. Math. Optim.} (2019), https://doi.org/10.1007/s00245-019-09555-4.

\bibitem{Trol}
F. Tr\"oltzsch,
Optimal Control of Partial Differential Equations. Theory, Methods and Applications,
{\it Grad. Stud. in Math.,} Vol. {\bf 112}, AMS, Providence, RI, 2010.

\bibitem{WLFC}
S.M. Wise, J.S. Lowengrub, H.B. Frieboes, V. Cristini,
Three-dimensional multispecies nonlinear tumor growth—I: model and numerical method. 
{\it J. Theor. Biol.} {\bf 253(3)} (2008), 524–-543.

\bibitem{WZZ}
X. Wu, G.J. van Zwieten and K.G. van der Zee, Stabilized second-order splitting
schemes for \CH~models with applications to diffuse-interface tumor-growth models, 
{\it Int. J. Numer. Meth. Biomed. Engng.} {\bf 30} (2014), 180-203.

\End{thebibliography}

\End{document}

\bye